\documentclass{amsart}
\usepackage{amscd,amssymb,subfigure,hyperref, epsfig}
\usepackage[arrow,matrix,graph,frame,poly,arc,tips]{xy}
\usepackage{amsmath}

\theoremstyle{plain}
\newtheorem{thm}[subsection]{Theorem}
\newtheorem{defn}[subsection]{Definition}

\newtheorem{prop}[subsection]{Proposition}
\newtheorem{cor}[subsection]{Corollary}

\theoremstyle{definition}

\newtheorem{exm}[subsection]{Example}

\numberwithin{equation}{section}

\newcommand{\HH}{{\mathcal H}}

\newcommand{\F}{{\mathbb F}}
\newcommand{\N}{{\mathbb N}} 
\newcommand{\Z}{{\mathbb Z}}

\newcommand{\R}{{\mathbb R}}

\newcommand{\Fqstar}{\F_q^{\,*}}

\renewcommand{\P}{{\mathbb P}}

\DeclareMathOperator{\rank}{rank}

\begin{document}

\title[Toric surface codes and Minkowski sums]%
{Toric surface codes and Minkowski sums}

\author[John Little]{John Little}
\address{Department of Mathematics and Computer Science,
College of the Holy Cross,  Worcester, MA 01610}
\email{\href{mailto:little@mathcs.holycross.edu}{little@mathcs.holycross.edu}}
\urladdr{\href{http://mathcs.holycross.edu/~little/homepage.html/}%
{http://mathcs.holycross.edu/\~{}little/homepage.html}}

\author[Hal Schenck]{Hal Schenck$^1$}
\address{Department of Mathematics,
Texas A\&M University, College Station, TX 77843}
\email{\href{mailto:schenck@math.tamu.edu}{schenck@math.tamu.edu}}
\urladdr{\href{http://www.math.tamu.edu/~schenck/}%
{http://www.math.tamu.edu/\~{}schenck}}

\thanks{$^1$Partially supported by NSF DMS 03-11142, NSA MDA
  904-03-1-0006, and ATP 010366-0103.}



%\date{May 24, 2006}

\begin{abstract}
Toric codes are evaluation codes obtained from an integral convex polytope
$P \subset \R^n$ and finite field $\F_q$.  They are, 
in a sense, a natural extension of Reed-Solomon codes, and have
been studied recently in \cite{dgv}, \cite{jh}, \cite{jh1}, 
and \cite{j}. In this paper, we obtain upper and lower bounds
on the minimum distance of a toric code constructed
from a polygon $P \subset \R^2$ by examining {\it Minkowski 
sum} decompositions of subpolygons of $P$. 
Our results give a simple and unifying explanation of bounds 
in \cite{jh1} and empirical results in \cite{j}; they also apply to 
previously unknown cases.
\end{abstract}
\maketitle

\section{Introduction}\label{sec:one}

In \cite{jh}, J. Hansen introduced the notion of a toric
surface code. Let $P \subset \R^2$ be an integral convex
polygon, and $\F_q$ a finite field such that 
after translation $P \cap \Z^2$ is properly contained in
the square $[0,q-2] \times [0,q-2]$
with sides of length $q-1$, which we denote $\Box_{q-1}$. 
Then a code is obtained
by evaluating monomials with exponent vector in $P \cap \Z^2$ 
at some subset (usually all) of the points of $(\Fqstar)^2$. 
We formalize this:
\begin{defn}
Let $\F_q$ be a finite field with primitive element
$\xi$. For $0 \le i,j \le q-2$ let $P_{ij} = (\xi^i,\xi^j)$ 
in $ (\Fqstar)^2$. For each $m = (m_1,m_2) \in P \cap \Z^2$, 
let $$e(m)(P_{ij}) = (\xi^i)^{m_1} (\xi^j)^{m_2}.$$ The toric
code $C_P(\F_q)$ over the field
$\F_q$ associated to $P$ is the linear code of 
block length $n = (q-1)^2$ spanned by the vectors in 
$\{(e(m)(P_{ij}))_{0\le i,j \le q-2} \colon m \in P \cap \Z^2\}.$ 
If the field is clear from the context, we will often omit it
in the notation and simply write $C_P$.
\end{defn}
The properties of these codes are closely tied to the geometry 
of the toric surface $X_P$ associated to the 
normal fan $\Delta_P$ of the polygon $P$. For example,
intersection theory on $X_P$ can be used to derive 
information about the minimum distance of toric codes.
The monomials $e(m)$ which are evaluated to produce the
generating codewords correspond to the lattice points 
$P \cap \Z^2$ and can be interpreted as sections of a certain 
line bundle on $X_P$. In \cite{jh1}, J. Hansen studies several specific 
families of polygons, depicted in Figure 1 below (notice that some families are 
completely contained in others). The minimum distance for these codes 
is determined by exhibiting codewords of weight equal to a lower bound 
obtained from intersection theory.
\begin{center}
\epsfig{file=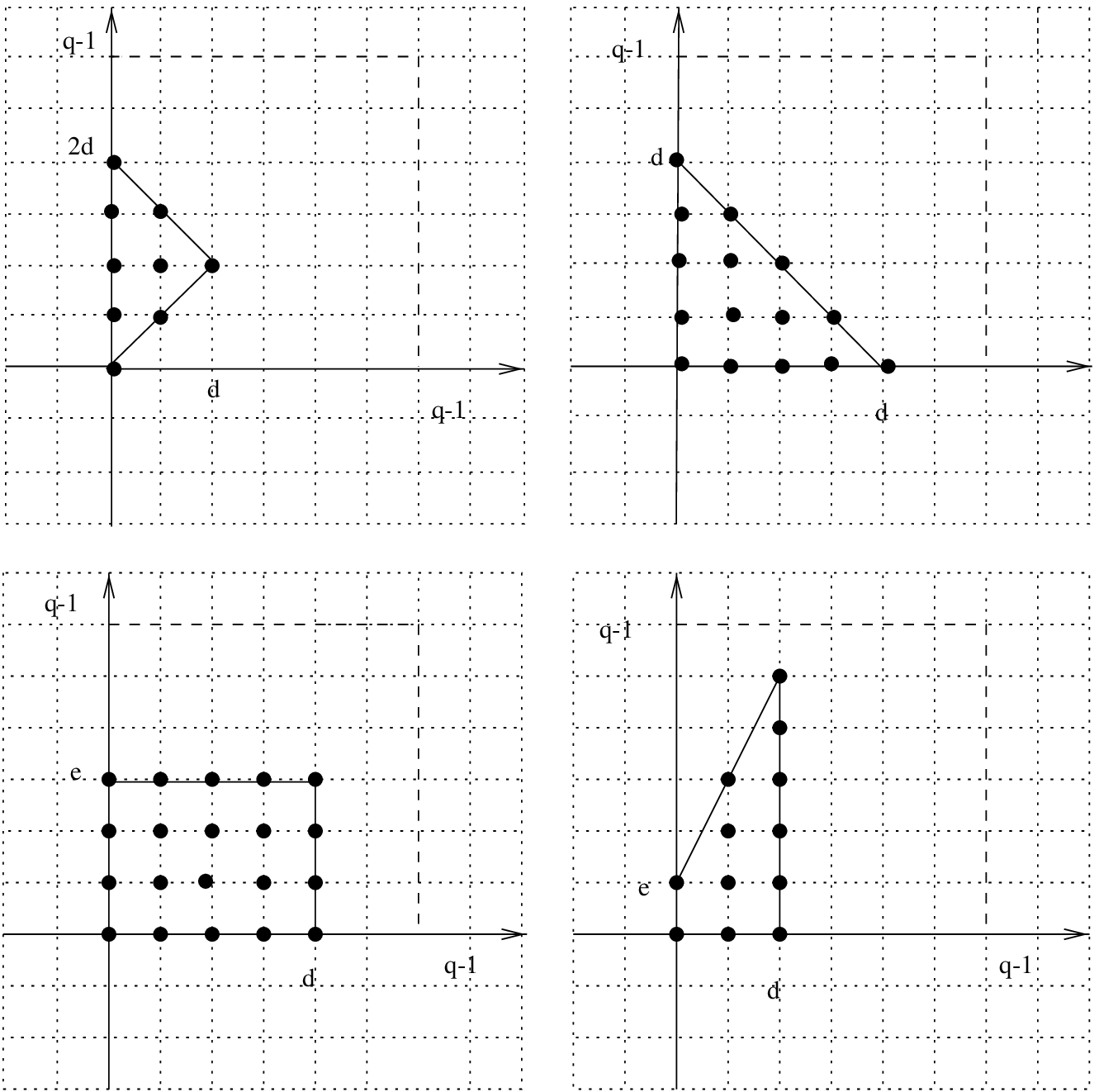,height=2.6in,width=3in}\\
Figure 1.
\end{center}
In this paper, we give upper and lower bounds on the minimum distance
for toric surface codes. Our formulas generalize the results of \cite{jh1},
and also provide theoretical explanations for the some of the values
tabulated in \cite{j}. Codewords
of small weight come from sections of the corresponding line
bundle that have many zeroes in $(\Fqstar)^2$.  A natural
way to try to obtain these is to consider sections that factor
into products of sections of related bundles
(we will call these {\it reducible} sections in the following).  
Such reducible sections come from polygons $P' \subseteq P$ that
decompose as {\it Minkowski sums} of other smaller polygons.

The definition of the Minkowski sum of polytopes will be 
reviewed in \S 2 below.  In Proposition ~\ref{upper}, we derive
an upper bound on the minimum distance of a toric surface code
when $P$ has a subpolygon that decomposes as a Minkowski sum
of other polygons.
We then apply these methods in \S 3 and \S4
to study the minimum distances of a number
of examples, including all toric surface codes from 
smooth toric surfaces $X$ with $\rank \mbox{Pic}(X) = 2$, or $3$.

In \S 5, we derive a statement complementary to 
the upper bound of Proposition ~\ref{upper}, giving a lower bound
on the minimum distance of toric codes constructed
from Minkowski-decomposable polygons.
The Hasse-Weil bound on the number of $\F_q$-rational points 
on a curve shows that for any given polygon $P$, there exists a 
lower bound on $q$ such that reducible sections of the corresponding
line bundle necessarily have more zeroes in $(\Fqstar)^2$ than irreducible 
sections.  For precise statements here, see Proposition~\ref{reduciblesections} 
and Corollary~\ref{ratsects} below.  This leads to our main theorem.

\begin{thm}
\label{mainth}
Let $\F_q$ be a finite field and let $P \subset \R^2$
be an integral convex polygon strictly contained in $\Box_{q-1}$.
Assume that $q$ is sufficiently large $($i.e. the bound (1) from Proposition 5.2 applies$)$.
Let $\ell$ be the largest positive integer such that there is some
$P' \subseteq P$ that decomposes as a Minkowski sum 
$P' = P_1 + P_2 + \cdots + P_\ell$
with nontrivial $P_i$. 
Then there exists some $P'\subseteq P$ of this form such that
$$d(C_P(\F_q)) \ge \sum_{i=1}^\ell d(C_{P_i}(\F_q)) - (\ell-1)(q-1)^2.$$
\end{thm}

\noindent
We then apply this result to some additional, less straightforward,
examples.
 
To relate our approach to other previous work, we note that
two very general methods for obtaining bounds on the 
minimum distance of codes on a higher dimensional variety $X$  
appear in recent work of S. Hansen \cite{sh}. The first method requires
finding the multipoint Seshadri constant for the line bundle whose 
sections are evaluated to obtain the code. The second method 
consists of covering the $\F_q$-rational points of $X$ with
curves and then counting points on these curves via 
inclusion-exclusion; of course, this depends on being able to find ``good'' 
curves on $X$. The methods we introduce here depend on finding sections 
which factor, so they relate to the second technique. 

The methods we use here make use of properties of toric surfaces in
an essential way. First, a key fact about {\it complete} toric varieties 
is that all the higher cohomology of a globally generated line bundle
vanishes. The lattice points in a polygon correspond to global
sections of such a line bundle, so Riemann-Roch provides a relation 
(see \S 5) between
lattice points and intersection theory. We also make use of the
Hasse-Weil bounds on the number of $\F_q$-rational points of a curve;
to apply the formula we need the arithmetic genus of an
irreducible section. The adjunction formula (\cite{f}, p. 91) 
gives the arithmetic genus in terms of polytopal data.

\section{Minkowski sums}\label{sec:two}

In this section, we give a brief discussion of the Minkowski 
sum operation, referring to Ziegler \cite{z} for more details. For
facts on toric varieties, our basic reference is Fulton \cite{f}.

\begin{defn}
Let $P$ and $Q$ be two subsets of $\R^n$. The Minkowski
sum is obtained by taking the pointwise sum of $P$ and $Q$:
$$P+Q = \{x+y\mbox{ } |\mbox{ } x \in P, y \in Q \}.$$
\end{defn}
We write conv to denote the convex hull of a set of points: 
the set of all convex combinations of the points.

\begin{exm}
Let $Q$ be the square $\mbox{conv}\{(0,0),(1,0),(0,1),(1,1)\}$ and 
let $P$ be the triangle $\mbox{conv}\{(0,0),(1,2),(2,1)\}$. Then
$$P+Q = \mbox{conv}\{(0,0),(1,0),(3,1),(3,2),(2,3),(1,3),(0,1)\}$$
as shown in Figure 2 below. 
\begin{center}
\epsfig{file=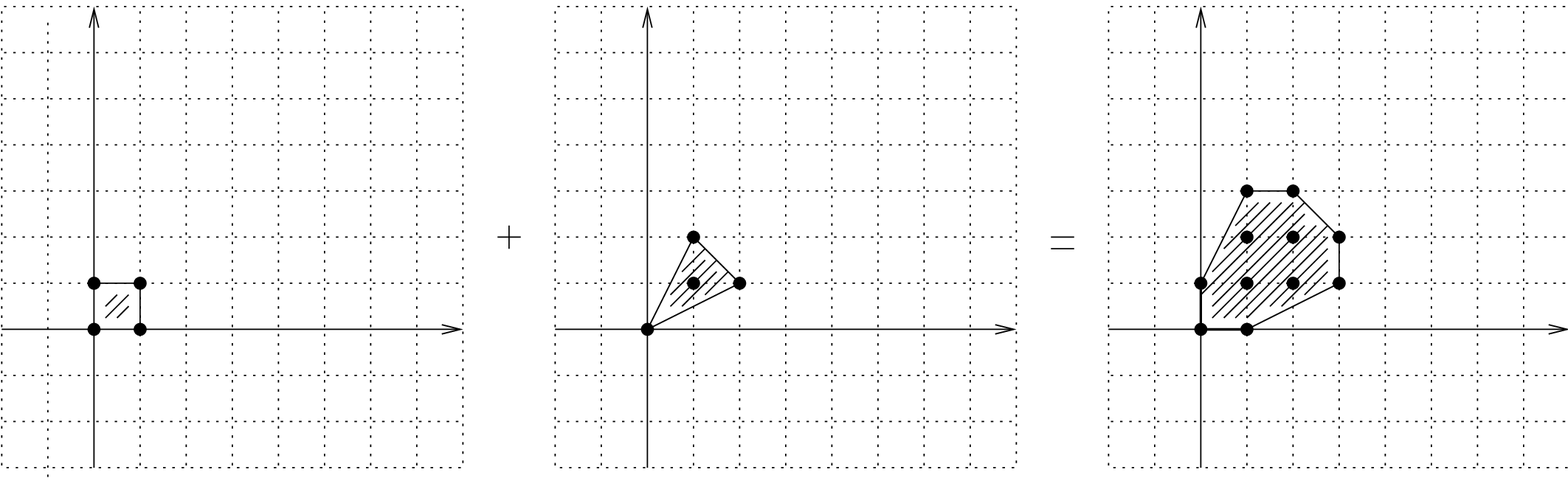,height=1.3in,width=4.5in}\\
Figure 2.
\end{center}
\end{exm}

If $f$ is a polynomial in two variables:
$$f(x,y) = \sum_{(a,b)\in \Z^2_{\ge 0}} c_{ab} x^a y^b,$$
then
$$NP(f) = \mbox{conv}\{(a,b) : c_{ab} \ne 0\}$$ 
is called
the {\it Newton polygon} of $f$.  It is a direct consequence
of the definition that if $f,g$ are two polynomials, then
$NP(fg) = NP(f) + NP(g)$, where the sum on the right is 
the Minkowski sum.   

Similarly, in the language of toric surfaces,
it is easy to see that if $P_1$ and $P_2$ are polygons,
then the normal fan $\Delta_{P_1+P_2}$ is the common
refinement of the fans $\Delta_{P_1}$ and $\Delta_{P_2}$. Thus, the
lattice points in $P_1+P_2$ correspond to
a basis of the global sections of a certain line bundle 
$\mathcal{O}(D)$ on the toric surface $X_{P_1+P_2}$, and the lattice
points in $P_1$ and $P_2$ 
correspond to bases of global sections for two other line bundles 
$\mathcal{O}(D_1)$ and $\mathcal{O}(D_2)$ on $X_{P_1+P_2}$ (see \cite{f}, p. 67).
If $D_1$ and $D_2$ are divisors on the toric surface $X$ corresponding
to polygons $P_1$ and $P_2$ with 
$s_1 \in H^0(\mathcal{O}(D_1))$ and
$s_2 \in H^0(\mathcal{O}(D_2))$ then
$$s_1s_2 \in H^0(\mathcal{O}(D_1)) \otimes H^0(\mathcal{O}(D_2)) \subseteq
H^0(\mathcal{O}(D_1 + D_2)),$$ 
which corresponds to the Minkowski sum $P_1 + P_2$ (indeed, if 
the $D_i$ are globally generated, then $H^0(\mathcal{O}(D_1))
\otimes H^0(\mathcal{O}(D_2)) = H^0(\mathcal{O}(D_1+D_2))$, \cite{f},
p. 69). A good exercise for toric experts is to work this out for the 
previous example. 

A first observation concerning the connection between the minimum 
distance of $C_P$ and Minkowski sums is the following.

\begin{prop}
\label{upper}
Let $\sum_{i=1}^\ell P_i \subseteq P$, and let $X$ be the 
toric surface corresponding to $P$.  Let
$m_i$ be the maximum number of zeroes in $(\Fqstar)^2$ 
of a section of the line bundle on $X$ corresponding to $P_i$, and 
assume that there exist sections $s_i$ with sets of $m_i$ zeroes
that are pairwise disjoint in $(\Fqstar)^2$.
Then 
$$d(C_P) \le \sum_{i=1}^\ell d(C_{P_i}) - (\ell - 1)(q-1)^2.$$
\end{prop}

\begin{proof}  By the definition we have
$d(C_{P_i}) = (q-1)^2 - m_i$.  As noted above, $N(fg)= N(f) + N(g)$, 
so the product $s = s_1 s_2 \cdots s_\ell$ is a section of the line bundle
$\mathcal{O}(D)$ corresponding to $\sum_{i=1}^\ell P_i$.
Moreover $s$ has exactly $m = m_1 + \cdots + m_\ell$ zeroes 
in $(\Fqstar)^2$ by hypothesis.  There is a 
codeword of the toric code $C_P$ with weight 
$$w = (q-1)^2 - m = \sum_{i=1}^\ell d(C_{P_i}) - (\ell - 1)(q-1)^2,$$
obtained by evaluating $s$. Hence
$$d(C_P) \le \sum_{i=1}^\ell d(C_{P_i}) - (\ell - 1)(q-1)^2,$$
which is what we wanted to show.
\end{proof}

Of course, the proof of the proposition can be extended
to handle the case where pairs of the $s_i$ have common zeroes in 
$(\Fqstar)^2$.  But the resulting bounds on $d(C_P)$ will involve 
the inclusion-exclusion principle and are harder to state in that
generality.  This upper bound also extends immediately to $m$-dimensional
toric codes for all $m\ge 2$ (that is, toric codes constructed from 
polytopes $P \subset \R^m$).

\section{First results and examples}\label{sec:three}

In this section we will present several results
on minimum distances of toric codes via Minkowski sum
decompositions.  These cases can be handled
without using Theorem~\ref{mainth}, and hence involve no hypothesis on $q$
other than that needed to ensure $P \subset \Box_{q-1}$. 

\begin{prop} 
\label{segment}
Let $P = \mbox{\rm conv}\{(0,0),(a,0)\}$ be a line
segment $($a one-dimensional polygon$)$.  Then for all $q > a + 1$,
$d(C_P(\F_q)) = (q-1)^2 - a(q-1).$
\end{prop}

\begin{proof}
The corresponding codes $C_P$ are products of $q-1$ 
copies of a Reed-Solomon code and the formula for 
the minimum distance follows directly.  Note that $P$
is also a Minkowski sum of $a$ line segments of length 1.
\end{proof}

\begin{prop}
\label{triangle}
Let $P$ be the integral triangle $P = \mbox{\rm
  conv}\{(0,0),(a,0),(b,c)\}$. If $a, b, c \ge 0$ and $a \ge b+c$, then 
for all $q > a+1$ $($so $P \subset \Box_{q-1})$,  
$$d(C_P(\F_q)) = (q-1)^2 - a(q-1).$$
\end{prop}

\begin{center}
\epsfig{file=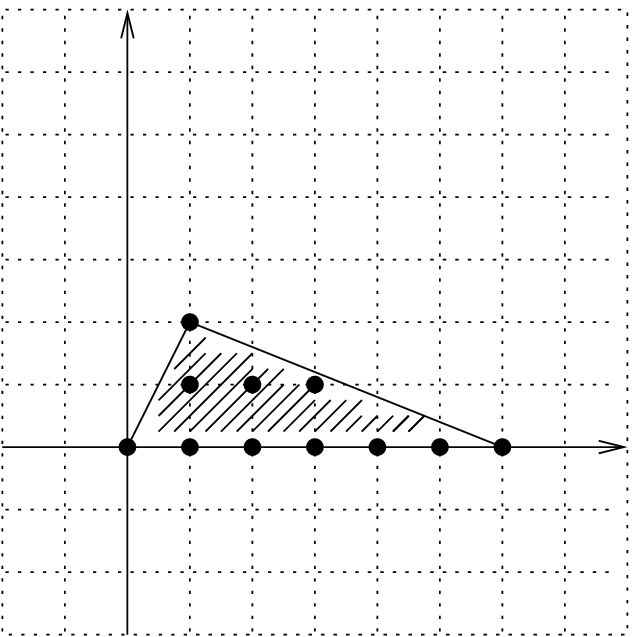,height=1.3in,width=1.5in}\\
Figure 3.
\end{center}

\begin{proof}
Note that $C_P$ may be viewed as a subcode of the code $C_{\Delta_a}$, where 
$$\Delta_a = \mbox{conv}\{(0,0),(a,0),(0,a)\}.$$
The toric surface corresponding to the triangle $\Delta_a$ 
is the $a$-tuple Veronese embedding of $\P^2$.
By a result of Serre (\cite{se}),
for all $q$, the curve of degree $a$ in $\P^2$ having 
the maximum possible number of 
$\F_q$-rational points is a reducible curve composed of $a$ 
concurrent lines.  When the point of intersection of the $a$ lines 
lies at infinity or on one of the coordinate axes in the
affine plane, then the corresponding curve
has $a(q-1)$ $\F_q$-rational points in $(\Fqstar)^2$.
Hence $d(C_P) \ge d(C_{\Delta_a}) = (q-1)^2 - a(q-1)$.  
Letting $P' = \mbox{conv}\{(0,0),(a,0)\}$, Proposition ~\ref{upper} shows
that $d(C_P) \le (q-1)^2 - a(q-1)$ as well.
\end{proof}

The code $C(\Delta_a)$ is also considered in \cite{jh1}, 
where the result $d(C_{\Delta_a}) = (q-1)^2 - a(q-1)$ is obtained in a 
different way.  

If $P'$ is any integral triangle obtained from $P$ by
a unimodular integer affine transformation (so $P$ and $P'$ are 
{\it lattice equivalent} polygons), then the same
formula applies to give $d(C_{P'})$.  This follows from 
the observation that if $P$ and $P'$ are lattice equivalent polygons, 
then $C_P$ and $C_{P'}$ are monomially equivalent
codes (\cite{ls}).  Propositions ~\ref{segment} and ~\ref{triangle} give a large collection of 
``building blocks'' to use in constructing other polygons.
We illustrate this by considering a standard class of toric surfaces
and toric codes studied in \cite{jh1}.

\begin{exm}
\label{Hirzebruch}
If $P = \mbox{conv}\{(0,0),(d,0),(0,e),(d,e+rd)\}$
for some $r \in \N$, then $P$ determines a {\it Hirzebruch surface},
denoted $\HH_r$. We assume $e + dr < q - 1$.
The polygon $P$ can be written as the Minkowski sum of
a line segment $L= \mbox{\rm conv}\{(0,0),(0,e)\}$ and a 
triangle $T= \mbox{\rm conv}\{(0,0),(d,0),(d,rd)\}$:
\vskip .075in
\begin{center}
\epsfig{file=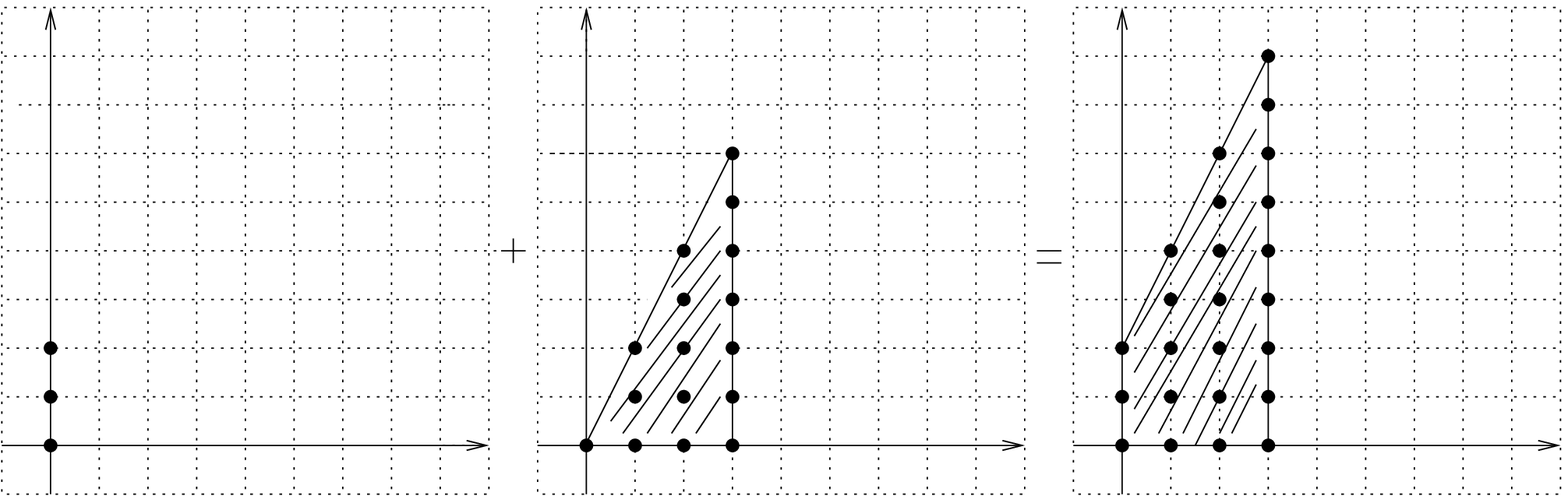,height=1.3in,width=4.5in}\\
Figure 4.
\end{center}
We now apply our results to this $P = T + L$ to determine the 
minimum distance of $d(C_P)$. The triangle $T$ is lattice equivalent 
to ${\rm conv} \{(0,0),(rd,0), (0,r) \}.$
By Proposition ~\ref{triangle}, for all $q$, $d(C_T) = (q-1)^2 - rd(q-1)$.
(The reducible sections of the line bundle corresponding to $T$
defined by $x^d\prod_{j=1}^{rd}(y - \alpha_j)$, $\alpha_j$ distinct in $\Fqstar$,
have exactly $rd(q - 1)$ zeroes in $(\Fqstar)^2$.
The $x^d$ corresponds to a trivial Minkowski summand and does not contribute
to the minimum distance.)  Similarly, Proposition ~\ref{segment}
shows $d(C_L) = (q-1)^2 - e(q-1)$. Thus, by Proposition ~\ref{upper} 
$$d(C_P) \le (q-1)^2 - e(q-1) + (q-1)^2 - rd(q-1)  - (q-1)^2 = (q-1)^2 - (rd + e)(q-1).$$
The polygon $P$ is a subset of a polygon lattice equivalent to 
the equilateral triangle $\Delta_{rd+e}$.   Hence $C_P$ is 
monomially equivalent to a subcode of $C_{\Delta_{rd + e}}$. 
It follows that the opposite inequality also holds, hence
$$d(C_P) = (q-1)^2 - (rd + e)(q-1).$$

Theorem 1.5 of \cite{jh1} gives $d(C_P)$ for the
codes from the Hirzebruch surfaces $\HH_r$ as the minimum of 
two terms.   Since the first term
given there is always larger than the second if $r > 0$, 
the minimum distance we obtain from the Minkowski sum decomposition agrees 
exactly with the value given in \cite{jh1}.
If $r = 0$, then the triangle $T$ reduces to a horizontal line segment,
and the Minkowski sum $T + L$ is a $d\times e$ rectangle.  The
corresponding toric code has minimum distance 
$$d(C_P) = (q-1)^2 - (d + e)(q-1) + de.$$
(see \cite{jh1}).  The minimum weight codewords
come from evaluating reducible sections
$$\prod_{i=1}^d (x - \alpha_i) \prod_{j=1}^e (y - \beta_j),$$
where the $\alpha_i$ are distinct and the  $\beta_j$ are distinct in $\Fqstar$.  
Note that this is one case where the factors have common zeroes, so Proposition
~\ref{upper} does not apply directly.
\end{exm}

For future reference, we note that by a result of Arkinstall (\cite{a}), 
the only lattice polygons with no interior lattice points are triangles
lattice equivalent to $\mbox{conv}\{(0,0),(p,0),(0,1)\}$
for some $p\ge 1$ or $\Delta_2$, or quadrilaterals with
two parallel sides.  Any such quadrilateral is 
lattice equivalent to one of the 
quadrilaterals defining a Hirzebruch surface with $d = 1$, or
to a $1\times e$ rectangle for some $e \ge 1$.  Hence by our discussion in Example
~\ref{Hirzebruch}, we know $d(C_P)$ for all toric codes from polygons 
$P$ with no interior lattice points.  

\section{Further examples: smooth surfaces with $\rank\,\mbox{Pic}(X) = 3$}\label{sec:four}

The Hirzebruch surfaces from Example ~\ref{Hirzebruch} satisfy $\rank\,\mbox{Pic}(\HH_r) = 2$
and, up to isomorphism, account for all smooth toric surfaces with this property.
In this section, we work out another extended family of examples and 
study the toric codes from the next most complicated
toric surfaces, those with $\rank\,\mbox{Pic}(X) = 3$.  
We will use some facts about toric surfaces, and refer to Section 2.5 
of \cite{f} for proofs. Recall that any smooth
complete toric surface $X$ may be obtained from $\mathbb{P}^2$ or
some $\HH_r$ by a succession of blow-ups at torus-fixed points.
The Picard number of such a surface is $n-2$, where $n$ is the 
number of 1-dimensional cones in the fan defining $X$. 
This description makes it reasonably straightforward to 
write down the fans for all smooth 
complete toric surfaces with $\rank\,\mbox{Pic}(X)=3$; either
we add a single ray to the fan of $\HH_r$ or a pair of rays to the
fan for $\mathbb{P}^2$, in such a way that for any two adjacent
rays, the determinant of the corresponding two by two matrix is $\pm 1$.

\begin{center}
\epsfig{file=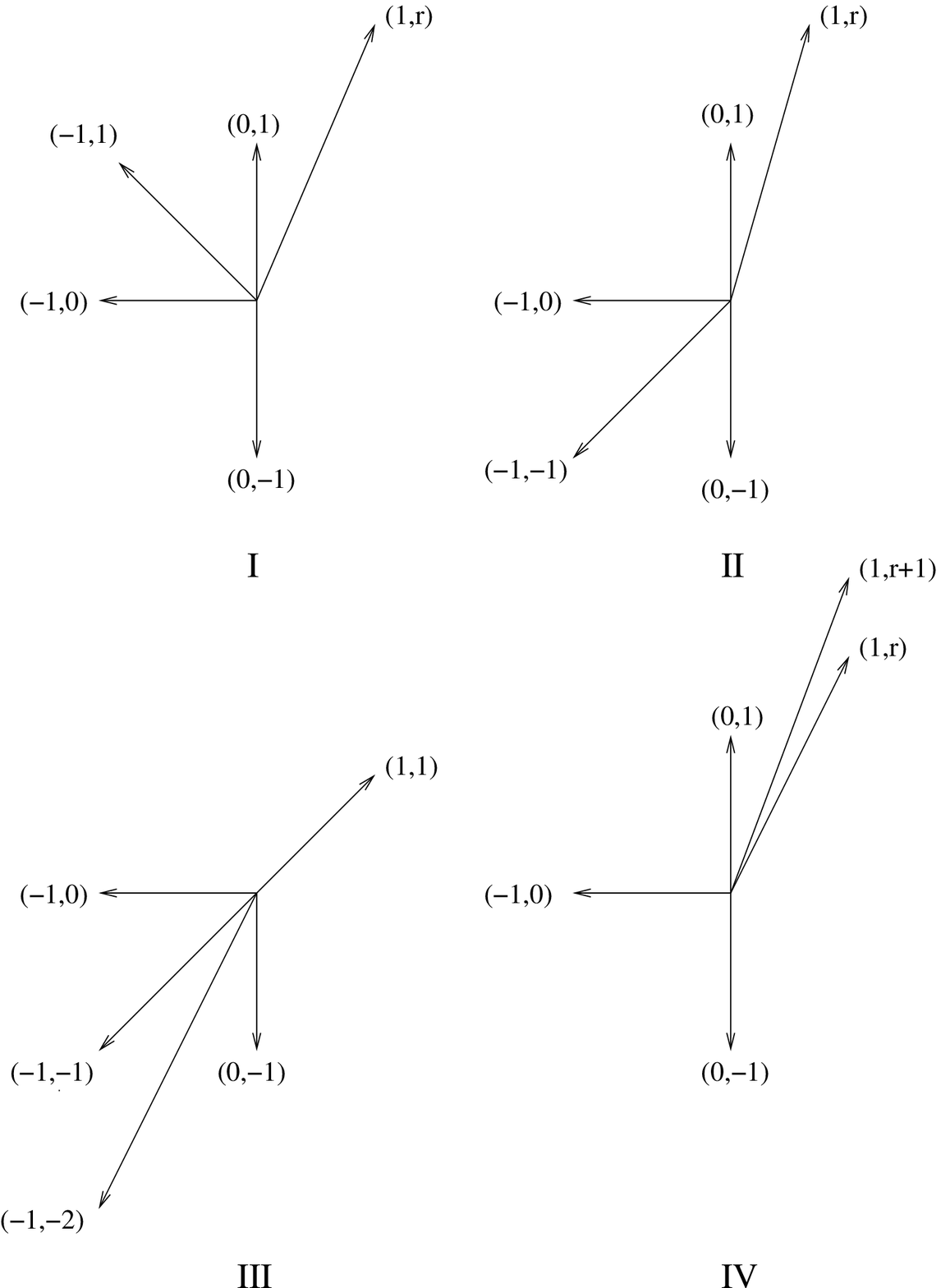,height=4.4in,width=3.3in}\\
Figure 5.
\end{center}

These fans are the outer normal fans of 
families of polygons. Polygons with these normal fans can ``scale'' 
in size, for example, the fan
with rays $\{(\pm 1,0),(0,\pm 1)\}$ is the normal fan for
any rectangle of the form $conv\{(0,0),(a,0),(0,b),(a,b)\}$.
In other words, the polytopes vary with parameters. 
We will see in a moment that these polygons all have Minkowski sum 
decompositions as sums of triangles and lines. 

For each fan, we want to determine the polygons whose
edges are normal to the given rays in the fan.  
For instance, in case I, polygons with this outer normal fan 
are obtained as the sets of solutions of inequalities as follows:
\begin{eqnarray*}
(1,r)\cdot (x,y) &\ge& \alpha\\
(0,1)\cdot (x,y) &\ge& \beta\\
(-1,1)\cdot (x,y) &\ge& \gamma\\
(-1,0)\cdot (x,y) &\ge& \delta\\
(0,-1)\cdot (x,y) &\ge& \varepsilon.
\end{eqnarray*}
for some $\alpha,\beta,\gamma,\delta, \varepsilon \ge 0$.
Taking $\delta = \varepsilon = 0$, $\gamma = a > 0$, $\beta = a + b$ with $b > 0$,
and $\alpha = r(a+b) + b + c$ with $c > 0$, we get a polygon as in Figure
6 below.

Now we are ready to examine Minkowski sum decompositions for
polygons corresponding to the fans in Figure 5.  For instance, in 
case I, we find that the pentagon 
$$P = \mbox{conv}\{(0,0),(r(a+b) + b + c,0),(b+c,a+b),
(b,a+b),(0,a)\}$$
can be decomposed as a Minkowski sum of the triangles
$$P_1 = \mbox{conv}\{(0,0),(ra,0),(0,a)\},\quad P_2 = \mbox{conv}\{(0,0),(b(r+1),0),(b,b)\}$$
and the line segment
$P_3 = \mbox{conv}\{(0,0),(c,0)\}$.
There are similar decompositions in each of the other cases
as well.
\vskip .2in
\begin{center}
\epsfig{file=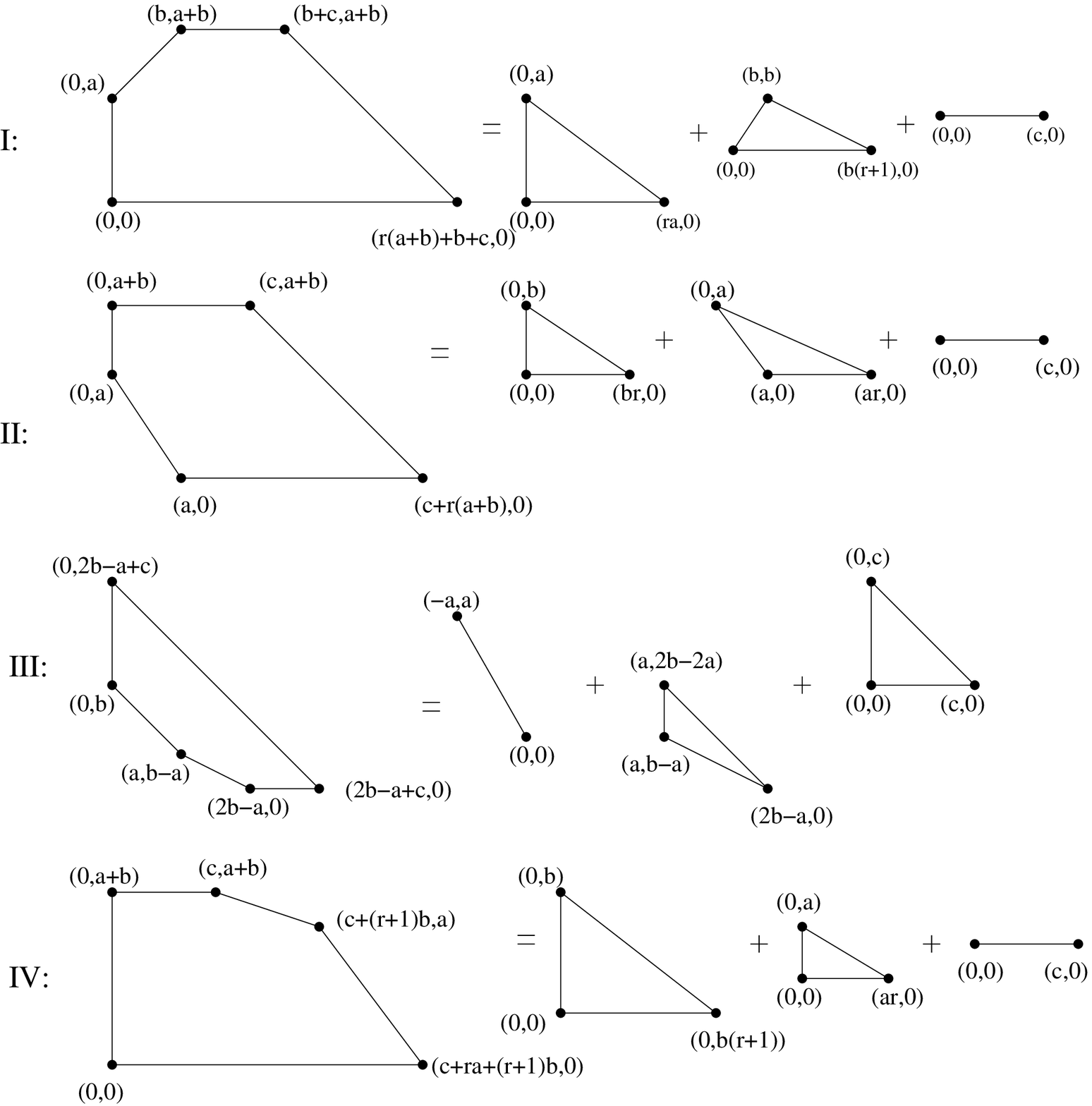,height=4.35in,width=4.6in}\\
Figure 6.
\end{center}

\begin{thm} Consider toric surface codes corresponding
to the families of polygons I,II,III,IV above, where $a,b,c,r \ge 1$
are integers and $q$ is sufficiently large so that the polygon
is contained in $\Box_{q-1}$.
\begin{enumerate}
\item In case I, for all such $q$, 
$$d(C_P) = (q-1)^2 - (r(a+b) + b + c)(q-1).$$
\item In case II, for all such $q$,
$$d(C_P) = (q-1)^2 - m(q-1),$$
where $m = \max\{a+b, c + (r-1)a + rb\}$.  
\item In case III, if $b > a$ as in Figure 7, then 
$$d(C_P) = (q-1)^2 - (2b+c-a)(q-1).$$
\item In case IV, for all such $q$,
$$d(C_P) = (q-1)^2 - (c + ra + (r+1)b)(q-1).$$
\end{enumerate}
\end{thm}

\begin{proof} We sketch how
the value in case I can be established using the methods presented in \S\S 2,3. 
We see first that the
stated value is an upper bound for $d(C_P)$ using the 
Minkowski sum decomposition given in Figure 6, Proposition
~\ref{upper}, Proposition ~\ref{segment} for the line segment and 
Proposition ~\ref{triangle} for the triangles.  Then, the fact that the given value for
$d(C_P)$ is the exact minimum
distance follows as in Example ~\ref{Hirzebruch}.  The polygon here
is contained in the equilateral triangle $\Delta_{r(a+b) + b + c}$.
Hence the minimum distance for $C_P$ is bounded
below by the minimum distance for the code
$C_{\Delta_{r(a+b) + b + c}}$.  
We leave it as an exercise for the reader to provide detailed 
proofs for the other parts. 
In each case, the minimum weight codewords come from
evaluation of reducible sections of the corresponding
line bundles.  For instance, 
the sections of $\mathcal{O}(D)$ with the maximal number of
$\F_q$-rational points in case III are given 
by $(y - \alpha_1x)\cdots (y - \alpha_{2b+c-a}x)$
with $\alpha_i \in \Fqstar$ distinct.  
\end{proof}

\section{Main theorem}\label{sec:five}

In this section we prove our main result, Theorem~\ref{mainth}.
The essential idea is to combine the Minkowski sum 
construction with the Hasse-Weil bounds on the number 
of $\F_q$-rational points of a curve:  If $Y$ is a 
smooth, absolutely irreducible curve over $\F_q$, then  
$$1+q-2g\sqrt{q} \le |Y(\F_q)|\le 1+q+2g\sqrt{q},$$
where $g$ is the genus of $Y$.  In \cite{ap}, the same inequalities
are demonstrated for absolutely irreducible, but possibly
singular curves, provided that $g$ is interpreted as the 
{\it arithmetic genus} of $Y$.  

The intuition behind our results is quite simple: 
from the Hasse-Weil bound, if $P$ is fixed and $q$ is 
sufficiently large, then sections which are reducible
must have more zeroes than irreducible sections.  

We will use the following notation.  For a polygon $P$, 
$v(P)$ will denote the area (2-dimensional volume) of $P$, 
$\#(P) = |P \cap \Z^2|$ will denote the
number of lattice points in $P$, $\partial(P)$ will
denote the number of 
lattice points in the boundary of $P$, and 
$I(P) = \#(P) - \partial(P)$ will denote the number of lattice 
points in the interior of $P$. Pick's theorem for
lattice polygons in $\R^2$ is the equality:
$$v(P) = \#(P) - \frac{1}{2} \partial(P) - 1.$$

Recall that we have seen that all polygons $P$ 
with $I(P) = 0$ correspond to toric surfaces
for which the minimum distance of $C_P$ is known
by results from \S 2, 3.  Hence, in the following
we will assume $I(P) > 0$.

In \S 2 we noted that if $\mathcal{O}(D_i), i \in \{1,\ldots, n \}$ are 
globally generated line bundles on a toric surface, then the global 
sections of $\mathcal{O}(\sum D_i)$ correspond to the Minkowski sum of
the polygons $P_i$ defined by $H^0(\mathcal{O}(D_i))$. 
Our starting data is a lattice polygon $P$, and to find 
reducible sections, our strategy is to work backwards: we look
for Minkowski sums $\sum_{i=1}^n P_i = P' \subseteq P$ with $n$ large.

In order to use algebraic geometry, we will first pass to a smooth 
surface. The toric surface $X_\Delta$ 
defined by the outer  normal fan $\Delta$ to $P$ need not be smooth.  
However, we can refine $\Delta$ to a fan $\Delta'$ such that 
$X_{\Delta'}$ is smooth, and the line bundle $\mathcal{O}(D)$ 
on $X_{\Delta'}$ corresponding to $P$ is generated by global sections 
(see \cite{f} p.90 or \cite{bs}). The numerical invariants 
$D^2$ and $DK$ discussed in the next paragraphs have simple
interpretations on the smooth surface $X_{\Delta'}$; most
importantly, they depend only on $P$. 

Finally, when we
deal with subpolygons $P_i$ of $P$, in order to make the same
set up work, we will refine the fan $\Delta'$ to include the 
outer normals to $P_i$, and then further subdivide the result (for
smoothness) to a fan $\Delta''$. The key point is that (\cite{f} p. 73) 
the $P_i$ correspond to globally generated line bundles on the smooth
surface $X_{\Delta''}$. So henceforth we will be working with 
globally generated line bundles on a smooth toric surface.

\begin{prop}
\label{algfacts}
Let $X$ be a smooth toric surface, and $K = K_X$ a canonical divisor.
Let $C$ be an irreducible curve on $X$ of arithmetic genus $g_C$ such
that the corresponding line bundle is globally generated, with 
$P$ the polytope corresponding to $H^0(\mathcal{O}(C))$. Then:
\begin{enumerate}
\item $g_C = \frac{C^2+CK}{2}+1.$

\item $h^0(\mathcal{O}(C))=\frac{C^2-CK}{2}+1.$

\item $g_C = 2v(P) + 2 - \#(P) = I(P).$
\end{enumerate}
\end{prop}
\begin{proof}
The first formula is simply adjunction, see \cite{h}, V.1.5 or
\cite{f}, p. 91. 
Since all the higher cohomology of a globally generated line bundle on 
a toric variety vanishes, and because a toric surface is rational, 
if $\mathcal{O}(C)$ is globally generated, then Riemann-Roch for
surfaces (\cite{h}, V.1.6) yields the second formula.
Adding the first two formulas shows that 
$h^0(\mathcal{O}(C)) + g_C = C^2+2.$ Since 
$h^0(\mathcal{O}(C)) = \#(P)$ and $C^2 = 2v(P)$ 
(see \cite{f}, p. 111), the last formula follows from Pick's theorem.
\end{proof}

One other fact that will be useful for us is that on a smooth 
toric surface $X$, the anticanonical divisor class $-K$  
is given by the sum of the divisors corresponding to the 
1-dimensional cones in the fan defining $X$ (\cite{f}, p. 85).
Now, 
$$(\F_q^*)^2 = X \setminus \bigcup\limits_{\tau \ne \{0\}} V(\tau),$$
where $V(\tau)$ is the closure of the torus orbit of the cone $\tau
\subseteq \Delta$, see \cite{f}, 3.1. In particular, a toric surface
decomposes as the union of a two dimensional torus with a finite set
of curves, which correspond exactly to the rays of $\Delta$.
Hence, the intersection number $-KC$ accounts for points on $C$ 
in the complement of the torus in $X$.

Our first result shows that if $q$ is sufficiently large,
then reducible sections with more irreducible components
necessarily have more zeroes in $(\Fqstar)^2$ than
sections with fewer irreducible components. In what 
follows, we write $V(s)$ for the zero locus of a section $s$. 
 
\begin{prop} 
\label{reduciblesections}
Let $P$ be a lattice polygon in $\R^2$ with $I(P) > 0$, 
and let $P' = \sum_{i=1}^m P_i'$ and $P'' = \sum_{k=1}^\ell P_k''$ 
(with $P_i'$ and $P_k''$ nontrivial) be two polygons contained in $P$.
Let $X$ be a smooth toric surface obtained by refining the normal 
fan $\Delta$ of $P$ as described above, so that $P'$ and $P''$
correspond to line bundles $\mathcal{O}(D')$ and $\mathcal{O}(D'')$ on
$X$. Let $s' = s_1's_2'\ldots s_m' \in H^0(\mathcal{O}(D'))$ and
$s'' = s_1''s_2''\ldots s_\ell'' \in H^0(\mathcal{O}(D''))$ be reducible sections
with $V(s_i')$ and $V(s_k'')$ irreducible.
If $m > \ell$ and 
$$q \ge (4I(P) + 3)^2,\leqno(1)$$
then
$$|V(s') \cap (\Fqstar)^2| > |V(s'') \cap (\Fqstar)^2|.$$
\end{prop}

\begin{proof}  
Let $D_i'$ be the divisor corresponding to $V(s_i')$,
and $D_k''$ be the divisor corresponding to $V(s_k'')$.
We write $g_i = g(D_i')$ and $g_k'' = g(D_k'')$.  
Our starting point is the observation that
$$|V(s') \cap (\Fqstar)^2| \ge \sum_{i=1}^m \Big((q+1) - 2g_i'\sqrt{q}\Big)
-\sum_{i<j} D_i'D_j' + D'K.$$
This follows because
$$|V(s') \cap (\Fqstar)^2| = \sum_{i=1}^m |V(s_i')|
-T -B,$$
where $T$ is the number of common intersection points of the curves
inside the torus $(\Fqstar)^2$ and $B$ is the number of 
points of $D'$ in the ``boundary'' $X \setminus (\Fqstar)^2$. 
Since the number of common
intersection points of $D_i'$ and $D_j'$ is the intersection number
$D_i'D_j'$, $T \le \sum_{i<j} D_i'D_j'$. As noted earlier, 
the number of points of $D'$ outside the torus is $-D'K$, so that 
$B \le -D'K$ (note that $D_i'D_j'$ and $-D'K$ do not distinguish
$\F_q$ rational points, so they may well overcount). 
Substituting the Hasse-Weil lower bound 
$|V(s_i')| \ge q+1 - 2g_i'\sqrt{q}$ gives the result.
Similarly, by the Hasse-Weil upper bound,
$$\sum_{k=1}^\ell \Big((q+1) + 2g_k''\sqrt{q}\Big) \ge |V(s'') \cap (\Fqstar)^2|.$$
Hence if $q$ satisfies 
$$(m-\ell)(q + 1) > 2\left(\sum_i g_i' + \sum_k g_k''\right)\sqrt{q} + \sum_{i<j} D_i'D_j' - D'K,\leqno(2)$$
then the conclusion of the proposition follows. Write 
$$\beta = \frac{1}{m-\ell}\left(\sum_i g_i' + \sum_k g_k''\right).$$
By Proposition \ref{algfacts}.3, $g_i' = I(P_i')$, so 
$$\sum_i g_i' = \sum_i I(P_i') \le I(P') \le I(P),$$ 
and similarly for $\sum_i g_i''$. Because $m-\ell \ge 1$, we see that
$$\beta \le \frac{2}{m - \ell} I(P) \le 2 I(P). \leqno(3)$$ 
The inequality (2) is quadratic in $\sqrt{q}$, so by the quadratic formula,
(2) will hold if
\begin{eqnarray*}
\sqrt{q} &>& \beta + \sqrt{\beta^2 + \sum_{i<j} D_i'D_j' - D'K + 1}\\
&\ge& \beta + \sqrt{\beta^2 + \frac{1}{m-\ell}\left(\sum_{i<j} D_i'D_j' - D'K\right) + 1}.
\end{eqnarray*}
Since $D' = \sum D_i'$, we have 
$$\sum_{i<j} D_i'D_j' = \frac{(D')^2 - \sum_i (D_i')^2}{2}.\leqno(4)$$
Now we apply (4) and Proposition \ref{algfacts}:
\begin{eqnarray*}
\sum_{i<j} D_i'D_j' - D'K + 1&=& \frac{(D')^2 - D'K}{2} - \frac{\sum_i
  \Big((D_i')^2 + D_i'K\Big)}{2} + 1 \\
&=& h^0(\mathcal{O}(D')) - \sum_i g_i' + m \\
&=& (\#(P') - \sum_i I(P_i')) + m  \\
&\le& 2 \#(P).
\end{eqnarray*}
The last step follows because $m \le \#(P)$ (each time we add in
a new Minkowski summand, we get at least one new lattice point in the
Minkowski sum), and $(\#(P') - \sum_i I(P_i')) \le \#(P') \le \#(P)$.
Now we use the theorem of P.R. Scott (\cite{s}):
$$ \#(P) \le 3 I(P) + 7$$
for a lattice polygon $P$ such that $I(P) > 0$. From the above, we see 
$$\sum_{i<j} D_i'D_j' - D'K  + 1 \le 6 I(P) + 14.$$
Hence if the lower bound (1) holds, since $I(P)>0$ we have
\begin{eqnarray*}
\sqrt{q} &\ge& 2I(P) + 2I(P) + 3 \\
         &=& 2I(P) + \sqrt{ 4 I(P)^2+ 12 I(P)+ 9}\\
         &>& 2I(P) + \sqrt{ 4 I(P)^2 + 6 I(P)+ 14}\\
         &\ge& \beta + \sqrt{\beta^2 + \sum_{i<j} D_i'D_j' - D'K + 1}
         \mbox{ by } (3),
\end{eqnarray*}
which is what we wanted to show.  
\end{proof}

A number of very crude estimates were
used to show that (1) implies the conclusion here.
Our lower bound on $q$ will rarely be sharp.  Much smaller lower bounds on 
$q$ can be obtained if we know more about possible
factorizations of sections of $\mathcal{O}(D)$.  For instance, we
have the following statement. 

\begin{cor}
\label{ratsects}
In the situation of Proposition ~\ref{reduciblesections}, 
suppose that $g_i' = I(P_i') = 0$ and $g_k'' = I(P_k'') = 0$ 
for all $i,k$.  Then the conclusion of Proposition
~\ref{reduciblesections} holds for all $q > \#(P) + m$.
\end{cor}

\begin{proof} In this case $\beta = 0$ in the proof of Proposition ~\ref{reduciblesections}.
\end{proof}

Theorem~\ref{mainth} follows almost immediately from Proposition ~\ref{reduciblesections}.

\begin{proof} (of Theorem~\ref{mainth})  Let $d = d(C_P)$.
Given $P$, the proposition shows that under the hypothesis (1)
on $q$, the number of zeroes of a section can always be increased by 
finding a reducible section in $H^0(\mathcal{O}(D))$
with more nontrivial factors, if there is one.  
Hence the sections with the largest number of zeroes in $(\Fqstar)^2$
must come from nontrivial factorizations with the largest possible number of 
factors.  Say $s = s_1 s_2 \cdots s_m$ is a nonzero section with the
maximum number of zeroes $(q-1)^2 - d$.  Then counting the number of zeroes, 
$$(q-1)^2 - d \le \sum_{i=1}^m m_i,$$
where $m_i$ is the number of zeroes of $s_i$.  We have 
$d(C_{P_i}) \le (q-1)^2 - m_i$ for each $i$.  Hence
$$ \sum_{i=1}^m m_i \le m (q-1)^2 - \sum_{i=1}^m d(C_{P_i}).$$
Rearranging the inequalities gives
$$d \ge \sum_{i=1}^m d(C_{P_i}) - (m - 1)(q-1)^2$$
as claimed.
\end{proof}

We have not tried to account for common zeroes of the $s_i$ in the proof
of the Theorem.  Moreover, in applying this statement, 
it is important to realize that
there may be several different subpolygons with the maximal
number of Minkowski summands.  The bound in Theorem~\ref{mainth}
is only guaranteed to hold for the one that minimizes
$$\sum_{i=1}^m d(C_{P_i}) - (m - 1)(q-1)^2.$$

\begin{exm}  Consider the polygon 
$$P = Q_1 + Q_2 := \mbox{conv}\{(0,0),(1,1),(2,1),(1,2)\} + \mbox{conv}\{(0,0),(1,0)\}.$$
Taking $P' = P$ and 
$$P'' = P_1 + P_2 := \mbox{conv}\{(1,1),(1,2)\} + \mbox{conv}\{(0,0),(1,0)\}\subset P$$
gives two different Minkowski-decomposable 
subpolygons of $P$ with the same number $m = 2$ 
of nontrivial summands.  However, since $I(Q_1) = 1$, the sections
having Newton polygon equal to $Q_1$ have arithmetic genus 1 and can have more 
zeroes in $(\Fqstar)^2$ than the rational curves corresponding
to the summands in $P''$.   So in applying Theorem~\ref{mainth}
to this example, we should use the decomposition $P = Q_1 + Q_2$
rather than $P'' = P_1 + P_2$.
In fact, we see this already for fields such as $\F_8$, 
where $q$ is much smaller than the bound from 
Proposition ~\ref{reduciblesections}.
Indeed, by a Magma computation using the routines from \cite{j},
$$d(C_P(\F_8)) = 33,$$
while
$\sum_{i=1}^2 d(C_{Q_i}(\F_8)) - (q-1)^2 = 33$,
and $\sum_{i=1}^2 d(C_{P_i}(\F_8)) - (q-1)^2 = 35$.
\end{exm}

Next we will show that our results shed some additional light on 
the good examples of toric surface codes tabulated in \cite{j}.

\begin{exm}
In Example 3.9 of \cite{j}, Joyner gives an example of a toric code
over $\F_8$ with $k = 11$ and $d = 28$.  
These parameters were better than any known code in Brouwer's tables 
\cite{b} at the time his article was written.  The convex hull of the 
integral points is a triangle $P = \mbox{conv}\{(0,0),(1,4),(4,1)\}$.  
Note that $P$ contains a translate of the triangle $\Delta_3$.  Applying 
Propositions ~\ref{upper} and ~\ref{triangle}, we obtain $d(C_P(\F_q)) \le (q-1)^2 - 3(q-1)$ 
for all $q$, so $d(C_P(\F_8)) \le 28$.  The lower bound 
$d(C_P(\F_q)) \ge (q-1)^2 - 3(q-1)$ also holds for $q$ sufficiently
large by Theorem~\ref{mainth}.  Joyner's
computations show that this bound on $d$ 
is also valid for $q = 8$, but our general statements
are not quite strong enough to prove this.
\end{exm}  
  
The following example gives an indication of
some additional interesting behavior that can 
occur for small $q$.

\begin{exm}
Consider the polygon 
$$P = \mbox{conv}\{(1,0),(2,0),(0,1),(1,2),(3,2),(3,3)\}.$$ 
\begin{center}
\epsfig{file=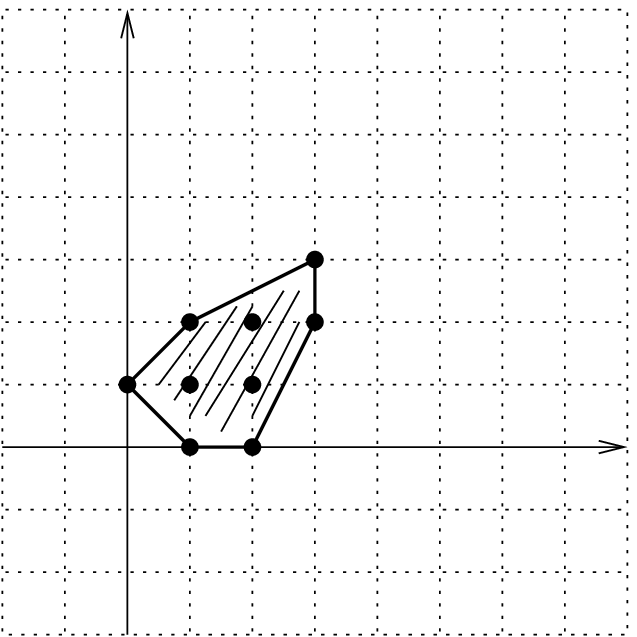,height=1.3in,width=1.5in}\\
Figure 7.
\end{center}
Note that $P \subset \Box_{q-1}$ for all $q \ge 5$.  
We see that $P$ contains the Minkowski-decomposable subpolygons
$P' = \mbox{conv}\{(1,0),(2,0),(1,2),(2,2)\}$ (a $1\times 2$ rectangle),
and $P'' = \mbox{conv}\{(1,0),(1,1),(3,2),(3,3)\}$
(a $2\times 1$ parallelogram).  $P'$ can be written as the 
Minkowski sum of two vertical line segments of length 1
and a horizontal line segment of length 1. Each $P_i$
gives $d(C_{P_i}) = (q-1)^2 - (q-1)$.  $P''$ has a similar
decomposition with three summands. 
There are no other Minkowski-decomposable subpolygons of $P$
with more than three Minkowski summands, and there are no
Minkowski summands with interior lattice points.  Hence we
have 
$$d(C_P(\F_q)) \ge (q-1)^2 - 3(q-1)$$
for $q > \#(P) + 3 = 12$ by Corollary ~\ref{ratsects}.

Both of these subpolygons give rise to 
reducible sections of the corresponding line bundles.
For instance from $P'$ we obtain reducible sections 
of the form $s = x(x-a)(y-b)(y-c)$. If $a,b,c\in \Fqstar$
and $b\ne c$, then $s$ has $3(q-1) - 2$ zeroes in $(\Fqstar)^2$.
Hence, by reasoning like that used in the proof of Proposition
~\ref{upper} (but in the case where the factors do have some common zeroes) we have 
$$d(C_P(\F_q)) \le (q-1)^2 - 3(q-1) + 2.$$

Computations using Magma show that 
\begin{eqnarray*}
d(C_P(\F_5)) = 6 & vs. & 4^2 - 3\cdot 4 + 2 = 6\\
d(C_P(\F_7)) = 20 & vs. & 6^2 - 3\cdot 6 + 2 = 20\\
d(C_P(\F_8)) = 28 & vs. & 7^2 - 3\cdot 7 + 2 = 30\\
d(C_P(\F_9))  = 42 & vs. & 8^2 - 3\cdot 8 + 2 = 42\\
d(C_P(\F_{11})) = 72 & vs. & 10^2 - 3 \cdot 10 + 2 = 72.
\end{eqnarray*}
The dimension is $k = \#(P) = 9$ in each case.  

The case $q = 8$ is the most interesting one here.  We may 
ask: Where does a section with $49 - 28 = 21$ zeroes in $(\F_8^{\,*})^2$
come from?  By examining the minimum weight codewords of this
code we find exactly 49 such words.  One of them comes, for
instance, from the evaluation of 
\begin{eqnarray*}
x + x^3y^3 + y^2 & \equiv & x(1 + x^2 y^3 + x^6 y^2) \bmod \langle x^7 - 1, y^7 - 1\rangle\\
                 & \equiv & x(1 + x^2 y^3 + (x^2 y^3)^3) \bmod \langle x^7 - 1, y^7 - 1\rangle
\end{eqnarray*}

Here $\langle x^7 - 1, y^7 - 1\rangle$ is the ideal of the
$\F_8$-rational points of the 2-dimensional torus.  So
$1 + x^2 y^3 + (x^2 y^3)^3$ has exactly the same zeroes in 
$(\F_8^{\,*})^2$ as $x + x^3y^3 + y^2$.  
Recall that $1 + u + u^3$ is one of the two irreducible polynomials
of degree 3 in $\F_2[u]$, hence 
$\F_8 \cong \F_2[u]/\langle 1 + u + u^3\rangle.$
Hence if $\beta$ is a root of $1 + u + u^3 = 0$ in $\F_8$, 
then 
$$1 + x^2 y^3 + (x^2 y^3)^3 = (x^2 y^3 - \beta)(x^2 y^3 - \beta^2)(x^2 y^3 - \beta^4)$$
and there are exactly $3\cdot 7 = 21$ points in $(\F_8^{\,*})^2$ 
where this is zero.  It is interesting to note that it is still a 
sort of reducibility that is producing a section with 
the largest number of zeroes here, even though the 
reducibility only appears when we look 
modulo the ideal $\langle x^7 - 1, y^7 - 1\rangle$.
We also note that these minimum weight codewords come from
curves with many rational points over the field $\F_8$
as in the construction used in \cite{bp}.  Similar phenomena
will occur for many other $P$ with $q$ small.
\end{exm}

\noindent{\bf Acknowledgments} This collaboration began while both
authors were members of MSRI during the commutative algebra program
in 2003. We also thank the Institute for
Scientific Computation at Texas A\&M for logistical support, and two
anonymous referees for careful readings and suggestions.

\bibliographystyle{amsalpha}

\begin{thebibliography}{10}

\bibitem{a} J. Arkinstall, {\em Minimal requirements for Minkowski's
    theorem in the plane}
            Bull. Austral. Math. Soc.  \textbf{22}  (1980), 259--283.

\bibitem{ap} Y.Aubry, M.Perret, {\em A Weil theorem for singular curves},
in Arithmetic, Geometry, and Coding Theory, R.Pellikaan, M.Perret,S.G.Vladut,
eds. de Gruyter, Berlin, 1996, 1-7.

\bibitem{b}  A.E. ~Brouwer, 
        {\em Bounds on linear codes},
             Handbook of Coding Theory, Elsevier (1998), 295-461, 
             and updates online at 
             {\tt http://www.win.tue.nl/\~{}aeb/voorlincod.html}.

\bibitem{bs}  T.~Beck, J. ~Schicho,
        {\em Sparse parametrization of plane curves},
             preprint, Radon Institute, 2005.

\bibitem{bp} P. ~Beelen, R. ~Pellikaan, {\em The Newton polygon of
plane curves with many rational points}, Des., Codes and Cryptography
\textbf{21} (2000), 41-67.

\bibitem{dgv} V. ~Diaz, C. ~Guevara, M. ~Vath,
        {\em Codes from n-Dimensional Polyhedra and n-Dimensional Cyclic Codes},
         Proceedings of SIMU summer institute, 2001.

\bibitem{f}  W. ~Fulton, 
        {\em Introduction to Toric Varieties},
        Princeton University Press, Princeton N.J., 1993.

\bibitem{jh}  J. ~Hansen, 
        {\em Toric surfaces and error-correcting codes}, in 
        Coding theory, cryptography and related areas (Guanajuato,
1998), 132--142, Springer, Berlin, 2000. 

\bibitem{jh1}  J. ~Hansen, 
        {\em Toric varieties Hirzebruch surfaces and error-correcting
codes},
        Appl. Algebra Engrg. Comm. Comput. \textbf{13} (2002),
289--300.

\bibitem{sh}  S. ~Hansen,
        {\em  Error-correcting codes from higher-dimensional
varieties},
        Finite Fields Appl. \textbf{7} (2001), 531--552.

\bibitem{h}  R. ~Hartshorne, {\em Algebraic Geometry}, Springer, 
New York, 1977.

\bibitem{j}  D. ~Joyner,
        {\em  Toric codes over finite fields},
        Appl. Algebra Engrg. Comm. Comput. \textbf{15} (2004),
63--79.
        
\bibitem{ls} J. ~Little, R. ~Schwarz, {\em On $m$-dimensional 
             toric codes}, {\tt arXiv:cs.IT/0506102}.

\bibitem{rab} S. ~Rabinowitz, {\em A census of convex lattice polygons 
with at most one interior lattice point},
             Ars Combin.  \textbf{28}  (1989), 83--96.

\bibitem{s} P.R. ~Scott, {\em On convex lattice polygons}, 
             Bull. Austral. Math. Soc. {\bf 15} (1976), 395-399.

\bibitem{se} J.P. ~Serre, {\em Letter \`a M. Tsfasman}, 
Journ\'ees Arithm\'etiques, 1989 (Luminy, 1989), Ast\'erisque {\bf 198-200} (1991), 351-353.
    
\bibitem{st} B. ~Sturmfels, 
        {\em Gr\"obner Bases and Convex Polytopes},
        AMS University Lectures Series, Vol. 8, 
        Am. Math. Soc., Providence, 1995. 

\bibitem{z} G. ~Ziegler, 
        {\em Lectures on Polytopes},
        Springer Verlag, Berlin, 1995. 

\end{thebibliography}

\end{document}